\documentclass[12pt]{amsart}
\newcommand{\Mo}{(M,\omega )}

\newcommand{\rar}{\rightarrow}

\usepackage{amssymb}
\usepackage{amsfonts}
\usepackage{amscd}
\newtheorem{thm}[subsection]{Theorem}
\newtheorem{lemma}[subsection]{Lemma}
\newtheorem{prop}[subsection]{Proposition}

\newtheorem{rem}[subsection]{Remark}

\newtheorem{problem}[subsection]{Problem}

\title{On simply connected $K$-contact non Sasakian manifolds}
\author{Bogus\l\/aw Hajduk}
\address{BH: University of Wroc\l\/aw, Plac Grunwaldzki 2/4, 50-384, Wroc\l\/aw, Poland}
\email{bhmath@interia.pl}
\author{Aleksy Tralle}
\address{AT: University of Warmia and Mazury, S\l\/oneczna 54, 10-710 Olsztyn, Poland}
\email{tralle@matman.uwm.edu.pl}
\begin{document}
\maketitle
\begin{abstract}{We solve the problem posed by Boyer and Galicki about the existence of simply connected $K$-contact manifolds with no Sasakian structure. We prove that such manifolds do exist using the method of fat bundles developed  in the framework of symplectic and contact geometry by Sternberg, Weinstein and Lerman.}

\end{abstract}

\section{Introduction}
In \cite{BG} Boyer and Galicki asked the following question (see Open Problem 7.4 on page 235 in this book).

\begin{problem}\label{probl:bg} Do there exist simply connected closed $K$-contact
manifolds with no Sasakian structure ?
\end{problem}

\noindent In this work, we answer this question positively.

\begin{thm}\label{thm:main} There exist simply connected
$K$-contact manifolds which do not carry any Sasakian structure.
\end{thm}

Let $(M,\eta)$ be a co-oriented contact manifold with a contact form
$\eta$. We say that $(M,\eta)$ is {\bf{\em $K$-contact}} if there is an
endomorphism $\Phi$ of $TM$ such that the following conditions are
satisfied:

\begin{enumerate}
\item $\Phi^2=-Id + \xi\otimes\eta,$ where $\xi$ is the Reeb
vector field of $\eta ;$

\item the contact form $\eta$ is compatible with $\Phi$ in the sense that
$$d\eta (\Phi X,\Phi Y)=d\eta (X,Y)$$ for all $X,Y$ and $d\eta (\Phi
X,X)>0$ for all nonzero $X\in \operatorname{Ker}\eta;$

\item the Reeb field of $\eta$ is a Killing vector field with respect to the
Riemannian metric defined by the formula
 $$g(X,Y)=d\eta (\Phi X,Y)+\eta
(X)\eta(Y).$$
\end{enumerate}

\noindent In other words, the endomorphism $\Phi$ defines  a complex
structure on $\operatorname{Ker}\eta$ compatible with $d\eta ,$
hence orthogonal with respect to  the metric $g = d\eta\circ
(\Phi\otimes Id).$ By definition, the Reeb field $\xi$ is orthogonal
to $\operatorname{Ker}\eta .$

\noindent For a contact manifold $(M,\eta )$ define  the {\bf{\em
metric cone}} or {\bf {\em the symplectization}}  as
$$\mathcal{C}(M)=(M\times\mathbb{R}^{>0},t^2\eta +dt^2).$$

\noindent Given a K-contact manifold $(M,\eta, \Phi ,g),$  the
almost complex structure $I$ on $\mathcal C(M)$ is defined by:

\begin{enumerate}

\item $I(X) = \Phi (X)$ on $\operatorname{Ker} \eta ;$

\item $I(\xi )= r\frac{\partial}{\partial r},
I(r\frac{\partial}{\partial r})=-\xi.$
\end{enumerate}

\noindent A K-contact manifold is called {\bf {\em Sasakian}},
if the almost complex
structure $I$ is integrable, hence defines a dilatation-invariant
complex structure on $\mathcal C(M),$ endowing $\mathcal{C}(M)$ with
a K\"ahler structure.

 Geometry of metric contact manifolds is important because
of their applications, for instance in the theory of Einstein
metrics \cite{BG}, \cite{BGM}.  $K$-contact manifolds have nice
topological properties (see \cite{BG}, Ch.7 and \cite{GNT}). For
example, they admit  Cohen-Macauley torus actions (this is an
analogue of the equivariant formality).  Most of the known examples of
$K$-contact manifolds are Sasakian, although examples of {\it
non-simply connected} $K$-contact manifolds with no Sasakian
structure are known \cite{BG}. This rises the following general
question: {\it given a contact manifold $M$, find conditions which
ensure that there exists a Sasakian metric compatible with the
contact structure.} This was posed by Ornea and Verbitsky in
\cite{OV}.

Although the main result of the paper belongs to the framework of
metric contact geometry, our methods come from symplectic and
contact geometry and are based on the notions of  symplectic and
contact fatness developed by Sternberg and Weinstein in the
symplectic setting \cite{W} and by Lerman in the contact case
\cite{L1}, \cite{L2}. Let  $G\rightarrow P\rightarrow B$ be a
principal bundle with a connection. Let $\theta$ and $\Theta$ be the
connection one-form and the curvature form of the connection,
respectively. Both forms have values in the Lie algebra
$\mathfrak{g}$ of the group $G$. Denote the pairing between
$\mathfrak{g}$ and   its dual $\mathfrak{g}^*$ by
$\langle\,,\rangle$.  By definition, a vector $u\in\mathfrak{g}^*$
is {\bf {\em  fat}}, if the two--form
$$
(X,Y)\rightarrow \langle\Theta(X,Y),u\rangle
$$
is nondegenerate for all {\it horizontal} vectors $X,Y$. Note that
if a connection admits at least one fat vector then it admits the
whole coadjoint orbit of fat vectors. In is important to notice that in this work we consider manifolds with contact forms rather than contact structures, but all the results from contact geometry which we use are valid in this particular situation.  Let $(M,\eta)$ be a contact
co-oriented manifold endowed with a contact action
 of a Lie group $G$. Define a {\bf{\em contact moment map}} by the formula
$$\mu_{\eta}: M\rightarrow \mathfrak{g}^*,\,\langle \mu_{\eta}(x),X\rangle=\eta_x(X^*_x)$$
for any $x\in M$ and any $X\in\frak{g}$. We denote by $X^*$ the
fundamental vector field on $M$ generated by $X\in\frak{g}$. Note
that the moment map depends on the contact form. The result below is
due to Lerman \cite{L2} and yields a construction of fiberwise
contact forms on the total space of the bundle associated to a
principal bundle with fat connection. \vskip6pt

\begin{thm}\label{thm:lerman-contact}

 Let there be given  a contact $G$-manifold $(F,\eta)$ with the contact moment map $\nu$.
Assume that
$$G\rightarrow P\rightarrow M$$
is a principal fiber bundle endowed with a connection such that the
image $\nu(F)\subset \frak{g}^*$ consists of fat vectors. Then there
exists a fiberwise contact structure on the total space of the
associated bundle

$$F\rightarrow P\times_GF\rightarrow M.$$

If the fiber $(F,\eta)$ is $K$-contact and $G$ preserves the
$K$-contact structure, then the total space of the associated bundle
is also $K$-contact.

\end{thm}
\noindent The second part of this theorem yields an explicit
construction of a fibered $K$-contact structure on a fiber bundle
and it is our tool to construct examples needed to prove Theorem
\ref{thm:main}. On the other hand, it is known that any closed
Sasakian manifold $M$ has Betti numbers $b_p(M)$ even for $p$ odd
and not exceeding $\frac12(dim\, M +1).$  We show, using Sullivan
models of fibrations \cite{FT}, that in our examples $b_3$ is odd,
hence they cannot be Sasakian.

Conceptually, contact manifolds are odd dimensional analogues of
symplectic manifolds, while manifolds with Sasakian structures are
treated as odd dimensional counterparts of K\"ahler manifolds. Thus,
Problem \ref{probl:bg} is an odd dimensional analogue of the topic
"symplectic vs. K\"ahler", which was quite important a decade ago
(\cite{CFM}, \cite{C}, \cite{IRTU}, \cite{RT}). Interestingly enough, like in the
symplectic case, the construction of non-simply connected
$K$-contact non-Sasakian manifold was found before the simply
connected one (compare \cite{TO}).

Finally, as a byproduct of our arguments and as a side remark we
answer positively a question posed by Alan Weinstein in 1980
\cite{W} about the existence of symplectically fat fiber bundle
whose total space is non-K\"ahler and simply connected. Note that
his question was asked in relation with a problem of constructing
simply connected non-K\"ahler symplectic manifold. As we have
already noted, the latter was solved by McDuff using the symplectic
blow-up construction \cite{McD}. Nevertheless, the question itself
remained open.

\section{Sullivan models}\label{sec:sullivan}
We use Sullivan models of fibrations as a tool of calculating
cohomology in low dimensions. In the sequel our notation follows
\cite{FT}. In this Section $\mathbb{K}$ denotes any field of zero
characteristic. We consider the
category of commutative graded differential algebras (or, in the
terminology of \cite{FT}, cochain algebras). If $(A,d)$ is a cochain
algebra with a grading $A=\oplus_pA^p$, the degree $p$ of $a\in A^p$
is denoted by $|a|$. \vskip6pt

Given a graded vector space $V,$ consider the algebra $\Lambda
V=S(V^{even})\otimes\Lambda (V^{odd})$, that is, $\Lambda V$ denotes
a free algebra which is a tensor product of a symmetric algebra over
the vector space $V^{even}$ of elements of even degrees, and an
exterior algebra over the vector space $V^{odd}$ of elements of odd
degrees.

We will use the following notation:
\begin{itemize}

\item by $\Lambda V^{\leq p}$ and $\Lambda V^{>p}$ are denoted the
subalgebras generated by elements of order $\leq p$ and of order
$>p,$ respectively;

\item if $v_i\in V$ is a generator, $\Lambda v$ denotes
the subalgebra generated by $v\in V$,

\item $\Lambda^pV=\langle v_1\cdot\cdot\cdot v_p\rangle$,
 $\Lambda^{\geq q}V=\oplus_{i\geq
q}\Lambda^iV$, $\Lambda^+V=\Lambda^{\geq 1}V.$
\end{itemize}
\noindent {\bf Definition}. A {\bf{\em Sullivan algebra}} is a
commutative graded differential algebra of the form $(\Lambda V,d)$,
where
\begin{itemize}
\item $V=\oplus_{p\geq 1}V^p$,
\item $V$ admits an increasing  filtration
$ V(0)\subset V(1)\subset...\subset V=\cup_{k=0}^{\infty}V(k)$
 with the property $d=0$ on $V(0)$,
$d: V(k)\rightarrow \Lambda V(k-1),\,k\geq 1$.
\end{itemize}

\vskip6pt \noindent {\bf Definition.} A Sullivan algebra $(\Lambda
V,d)$ is called {\bf{\em minimal}}, if
$$Im d\subset \Lambda^+ V\cdot\Lambda^+V.$$

\vskip6pt \noindent {\bf Definition}. A {\bf{\em Sullivan model}} of
a commutative graded differential algebra $(A,d_A)$ is a morphism
$$m: (\Lambda V,d)\rightarrow (A,d_A)$$
inducing an isomorphism $m^*: H^*(\Lambda V,d)\rightarrow
H^*(A,d_A)$. If $X$ is a CW-complex, there is a cochain algebra
$(A_{PL}(X),d_A)$ of polynomial differential forms. For a smooth
manifold $X$ we take a smooth triangulation of $X$ and as  the model
of $X$ the Sullivan model of $A_{PL}(X).$ If it is minimal, it is
called the {\it Sullivan minimal model of} $X$. We will need the
following known fact. If $(\frak{m}_X,d)=(\Lambda V,d)$ is the
minimal model of finite simply connected CW complex $X$, then

$$V\cong Hom_{\mathbb{Z}}(\pi_*(X),\mathbb{K}). \eqno(2)$$

\vskip6pt

\noindent {\bf Definition}. A {\bf{\em relative Sullivan algebra}}
is a graded commutative differential algebra of the form $(B\otimes
\Lambda V,d)$ such that
\begin{itemize}
\item $(B,d)=(B\otimes 1,d), H^0(B)=\mathbb{K},$
\item $1\otimes V=V=\oplus_{p\geq 1}V^p$,
\item $V=\cup_{k=0}^{\infty} V(k), V(0)\subset V(1)\subset...$,
\item $d: V(0)\rightarrow B, d(V(k))\rightarrow
B\otimes\Lambda V(k-1),\,k\geq 1$
\end{itemize}

\begin{lemma}\label{lemma:rel3} Consider a minimal Sullivan
algebra $(\Lambda V,d)$ such that $V^1=\{0\}$. Let $\omega$ be a
cocycle of degree 4 such that$[\omega]\in H^4(\Lambda V,d)$ is a
non-zero cohomology class. Define  the relative Sullivan algebra
$(\Lambda V\otimes\Lambda y,d)$ with generator $y$ of degree 3 by
$dy=\omega.$ Then
$$H^3(\Lambda V\otimes\Lambda y,d)\cong H^3(\Lambda V,d).$$
\end{lemma}

\noindent{\bf{\em Proof.}} Since $V^1=\{0\}$ and the degree of $y$
is $3$, one can write
$$dy=\sum v_i\cdot v_j+u,\,|v_i|=2,\,u\in  V^4.$$

If $z\in Z^3(\Lambda V\otimes\Lambda
y)$, one can write
$$z=x+\alpha y,\,\alpha\in\mathbb{K},\,x\in V^3.$$
Since $z$ is a cocycle, we have
$$0=dz=dw+\alpha(\sum v_iv_j+u).$$
Since $\Lambda V\otimes\Lambda y$ is a free algebra and $u$ is
indecomposable of degree $4$, thus, by minimality, either $u=0$, or
$\alpha=0$. Thus $z\in V^3.$ This means that  $$Z^3(\Lambda
V\otimes\Lambda y)=Z^3(\Lambda V),$$ and the proof follows.

\hfill$\square$

\begin{rem} {\rm One can consult \cite{TO} for a simple description
of a method of calculating minimal models of free cochain algebras.}
\end{rem}
\noindent Relative Sullivan algebras are models of fibrations. Let
$p: X\rightarrow Y$ be a Serre fibration with the homotopy fiber
$F$. Choose Sullivan models

$$m_Y: (\Lambda V_Y,d)\rightarrow (A_{PL}(Y), d_Y),\,\bar m: (\Lambda V,\bar d)\rightarrow A_{PL}(F).$$
There is a commutative diagram of cochain algebra morphisms

$$
\CD
A_{PL}(Y) @>>> A_{PL}(X) @>>> A_{PL}(F)\\
@A{m_Y}AA @A{m}AA @A{\bar m}AA\\
(\Lambda V_Y,d) @>>> (\Lambda V_Y\otimes\Lambda V,d) @>>> (\Lambda
V,\bar d)
\endCD
$$
in which $m_Y,m,\bar m$ are all Sullivan models (see Proposition
15.5 and 15.6 in \cite{FT}).

\section{Betti numbers of Sasakian manifolds}\label{betti}

In what follows all cochain algebras and cohomologies are assumed
to be defined over $\mathbb{K}=\mathbb{R}$. We will need the
following property of Sasakian
manifolds.

\begin{thm}\label{thm:oddbetti} If $M$ is closed Sasakian manifold of dimension
$2n+1,$  then for any odd $p\leq n+1$ the Betti numbers $b_p$ are
even.
\end{thm}

This was proved by several authors, cf. \cite{BG,F,T} and Theorem
7.4.11 in \cite{BG}. For  completeness as well as better
presentation of our methods we will explain a short proof.

Consider a smooth principal $S^1$-orbibundle $\pi : E\rar B$ over a
compact K\"ahler orbifold $B$ with smooth total space $E$ , i.e.,
$E$ is a smooth manifold, the circle acts smoothly on $E$ with all
isotropy groups finite and $\pi :E\rightarrow B$ is the natural
projection onto the orbit space of the action.  For the definition
and relevant properties of orbibundles see \cite{BG}, Chapter 4.

  Any compact K\"ahler orbifold satisfies hard Lefschetz
condition, i.e., it has a Lefschetz class $v\in H^2(B).$ By this we
mean that for $p=0,1,...n-1,\ dim\, B=2n,$ the linear maps
$$L_v^k: H^{n-p}(B)\rightarrow H^{n+p}(B):x\mapsto v^p\cup x $$
are isomorphisms (cohomology groups with real coefficients). There
is also the Gysin sequence for a sphere orbibundle $p:E\rightarrow
B.$ Both these properties are the consequences of the fact that the
category of K\"ahler orbifolds and orbibundles are rationally
equivalent to the category of K\"ahler manifolds and bundles.
Explicit proofs, using calculations of basic cohomology of the
foliation of $E$ by fibers, can be found in \cite{EK,NR,KT,WZ}.

\begin{lemma}\label{lemma:orbibun} If $p\leq n+1$ is odd, then the
Betti numbers $b_p(E)$ are even.
\end{lemma}

Proof. It follows from the Lefschetz property and duality that  for
a compact K\"ahler orbifold $B$ the Betti numbers $b_p(B)$ are even
for $p$ odd. The Gysin sequence for a principal circle orbibundle is
the long exact sequence

$$ ...\rightarrow H^{p-2}(B) \stackrel{L_v^1}{\rightarrow}
H^p(B)\stackrel{\pi^*}{\rightarrow} H^p(E)\rightarrow H^{p-1}(B)
\stackrel{L_v^1}{\rightarrow} H^{p+1}(B)\rightarrow ...$$

By assumption $L_v^i=L_v^{i-1}L_v^1$ is an isomorphism, thus for
$i\leq n$ the linear map $L_v^1:H^i(B)\rar H^{i+2}(B)$ is a
monomorphism. The Gysin sequence gives that $\pi^*$ is onto and
$H^p(E)=H^p(B){/}H^{p-2}(B)$ is even dimensional for any $p\leq
n+1.$ \hfill $\square$

Theorem \ref{oddbetti} now follows from the fact that any Sasakian
manifold is the total space of a circle orbibundle over a K\"ahler
orbifold. This can be found in \cite{OV}, where the authors  show
that if $M$ is Sasakian, than it admits also a quasi-regular
Sasakian structure. They  essentially repeated the argument of
Rukimbira \cite{R}, who noticed this for K-contact manifolds.  For a
quasi regular Sasakian manifold $M$ it had been known (and is rather
obvious) that $M$ admits a locally free circle action such that the
orbit space is a K\"ahler orbifold. Warning: in the statement of
Theorem 1.11 of \cite{OV} the word "fibration" should be understood
as "orbibundle". \hfill${\square}$

\section{Contact fat circle bundles}\label{sec:contact}

Let $(M,\omega)$ be a symplectic manifold with integral symplectic
form. Consider the principal circle bundle $\pi : P\rightarrow M.$
Note that in the case of circle bundles the curvature has real
values, if we identify the Lie algebra of $S^1$ with the reals. This
implies that if a principal $S^1-$bundle is fat, then only the zero
vector is non-fat. Assume that the bundle is determined by the
cohomology class $[\omega]\in H^2(M,\mathbb{Z})$. Fibrations of this
kind were first considered in \cite{BW} and are called {\bf{\em the
Boothby-Wang fibrations}}. By \cite{K},  $P$ carries an invariant
contact form $\theta .$ Namely, there exists an invariant connection
on $P$ with the curvature 2-form equal to $\pi^*\omega .$  Since the
Reeb vector field of such form is the infinitesimal generator of the
circle action on $P,$ the moment map is constant and nonzero.
Moreover, by \cite{W} principal circle bundles are fat if and only
if they are the Boothby-Wang fibrations.

Consider now  a $S^1-$contact manifold and the bundle associated to
$\pi .$ The following follows from \ref{thm:lerman-contact}, cf.
also \cite{BG}, ch.7.

\begin{thm} \label{thm:Boothby-Wang} Let there be given a Boothby-Wang fibration
$$S^1\rightarrow P\rightarrow M.$$
Assume that $(F,\eta, S^1)$ is a contact manifold endowed with an
$S^1$-action preserving  $\eta$ and the moment map has only non-zero
values. Then the associated bundle
$$F\rightarrow P\times_{S^1}F\rightarrow M$$
admits a fiberwise contact form. If $(F,\eta)$ is $K$-contact, the
same is valid for the fiberwise contact structure on
$P\times_{S^1}F$.
\end{thm}
\hfill$\square$

\section{A construction of simply connected $K$-contact non
Sasakian manifold}\label{sec:construct}

\subsection{Proof of Theorem \ref{thm:main}}
\begin{prop}\label{prop:construction} Let $(X,\omega )$ be a closed simply
connected symplectic manifold  and let $L\rightarrow X$ be the
complex line bundle corresponding to the cohomology class
$[\omega]\in H^2(X).$ Consider the Whitney sum $L\oplus L\rightarrow X$ and
the unit sphere bundle
$$S^3\rightarrow M \rightarrow X.$$
 The total space $M$ has a relative
Sullivan model of the form
$$(\frak{m}_X\otimes\Lambda y,d), dy=z\not=0,$$

\noindent where  $\frak m_X$ is a minimal model of $X,$
$z\in\frak{m}_X,$ $|z|=4,|y|=3$ and $z$ represents $[\omega ]^2.$
\end{prop}

\noindent{\bf{\em Proof.}} The Whitney sum $L\oplus L\rightarrow X$
is the rank $2$ complex vector bundle. From the multiplication
formula for the total Chern class of the Whitney sum of bundles one
obtains $c_2(L\oplus L)=c_1(L)^2=[\omega]^2\not=0$. This enables us
to apply the argument of Example 4 on page 202 of \cite{FT}. Namely,
any spherical fibration $p: Z\rightarrow X$ has the model
$$(A_{PL}(Z),\hat d)\simeq (A_{PL}(X)\otimes\Lambda y,d), dy=z\in A_{PL}(X).$$
If, moreover, the spherical fibration arises as a unit sphere bundle
of a vector bundle $\xi: E\rightarrow X$ of even real rank $k+1$,
then the cohomology class $[z]$ is the Euler class of $E$. Since the
given $S^3$-sphere bundle is the sphere bundle of $L\oplus L$, and
the latter has the Euler class equal to the second Chern class
$[\omega]^2\not=0$, we get the model of $M$ of the form

$$(A_{PL}(X)\otimes\Lambda y,d),\,dy=z\not=0, [z]=[\omega]^2.$$

\noindent By the theory described in Section \ref{sec:sullivan},
the above cochain algebra has the required  Sullivan model (of $M$).

\hfill$\square$

\begin{prop}\label{prop:example} Let $(X,\omega )$ be any compact simply
connected symplectic manifold such that $b_3(X)$ is odd. If
$$S^1\rightarrow P\rightarrow X$$
is a Boothby-Wang fibration with the Euler class equal to  $[\omega
],$  then the total space of the fiber bundle
$$S^3\rightarrow P\times_{S^1}S^3\rightarrow X$$
associated to $P\rightarrow X$ by the Hopf action of $S^1$ on $S^3$
admits a $K$-contact structure, but no Sasakian structure.
\end{prop}

\noindent{\bf{\em Proof.}} Boothby-Wang fibrations are fat and for
the contact form given by the curvature of an invariant connection
the moment map has nonzero values. Whence if $F$ is the total space
of such fibration, then the assumptions on $(F,\eta )$ in
\ref{thm:Boothby-Wang} are fulfilled. The sphere $S^3$ is the
Boothby - Wang fibration over $S^2,$ thus it is $K$-contact and the
Hopf $S^1$-action preserves the $K$-contact structure. This is the
special case of \ref{thm:Boothby-Wang}, with $F=S^1.$  However, it
is not difficult to give direct calculations (see below), and one
can also use  descriptions of $K$-contact manifolds in \cite{BG},
Th. 6.1.26, or in \cite{BMS}.  By Theorem \ref{thm:Boothby-Wang},
the total space $P\times_{S^1}S^3$ of the associated bundle admits a
$K$-contact structure as well. We claim that there exists a
fibration of this type which cannot be Sasakian. To prove the
latter, let us make the following observations. Denote by
$L\rightarrow X$ the complex linear bundle given by $P\rightarrow
X.$  Consider the Whitney sum $L\oplus L$ of $L$ as in Proposition
\ref{prop:construction}. This bundle  is the associated  to the
principal bundle $S^1\rightarrow P\rightarrow X$ if the circle acts
diagonally on $\mathbb C \oplus \mathbb C .$ If we pass to the
sphere bundle of $L\oplus L,$ then we see that this bundle is
associated to $P\rightarrow X$ by the Hopf circle action on $S^3$
and its total space is $M=P\times_{S^1}S^3.$ The Euler class of the
$S^3-$bundle is equal to $c_2(L\oplus L)=c_1(L)^2=[\omega ]^2\neq
0.$ This implies that $M$ has a relative Sullivan model given by
Proposition \ref{prop:construction}. Since $X$ is simply connected,
$(\frak{m}_X,d)=(\Lambda V,d)$ and, by formula $(2),$ we have
$V^1=\{0\} .$ Applying Lemma \ref{lemma:rel3} we get $b_3(M)=b_3(X)$
odd. This contradicts  Theorem \ref{thm:oddbetti} and the proof is
complete.

\hfill$\square$

Theorem \ref{thm:main} can be extended to get a whole class of
examples as follows.

\begin{thm} Let $(X,\omega )$ be a closed symplectic manifold and
$\xi : P \rightarrow X$ be the principal circle bundle with the
Euler class $[\omega ].$ Consider a unitary representation of $S^1$
in $\mathbb C^{n+1}$ such that all weights are positive and the
associated sphere bundle $E=P\times_{S^1}S^{2n+1}\rightarrow X.$
Then for all $n\geq 0$ the total space $E$ is K-contact. If
additionally $\pi_1(X)=1$ and $b_3(X)$ is odd, then $E$ admits no
Sasakian structure.
\end{thm}

Proof. On $S^{2n+1}\subset \mathbb C^{n+1}\cong \mathbb R^{2n+2}$ we
consider the standard K-contact form $\eta_{st} =
\iota_{\frac{\partial}{\partial r}}\omega_{st},$ where
$\omega_{st}=\sum_{i=1}^{n+1} dx_idy_i$ is the standard symplectic
form and $\frac{\partial}{\partial r}$ is the radial vector field.
The form $\eta_{st}$ is $U(n+1)-$invariant with the Reeb field
$$R=\sum_{i=1}^{n+1}(y_i\frac{\partial}{\partial x_i}-
x_ifrac{\partial}{\partial y_i}.$$ If $J$ is the standard complex
structure of $\mathbb C^{n+1},$ then $J(\frac{\partial}{\partial
r})=R$ and $J$ preserves $ker\, \eta_{st}.$ In particular, the
Riemannian metric $d\eta_{st}\circ(J\otimes Id)$ is the standard
round metric of $S^1,$ the restriction of the Euclidean metric of
$\mathbb R^{2n+2}.$ Moreover, $R$ generates the action of he circle
on $S^{2n+1}$ obtained from the diagonal representation in $\mathbb
C^{n+1}$ (with all weights equal to 1).

Consider now a unitary representation in $\mathbb C^{n+1}$ with
positive weights $w_1,...,w_{n+1}.$  We claim that the associated
bundle $S^{2n+1}\rightarrow E\rightarrow X$  satisfies assumptions
of Lerman's theorem. Since $\xi$ is fat and $\eta_{st}$ is
$U(n+1)-$invariant, we need only to check that the image of the
moment map does not contain zero. As $R$ is orthogonal to $ker\,
\eta_{st},$ this is equivalent to the non-vanishing of the scalar
product  $<R,V>,$ where $V$ generates the $S^1-$action on
$S^{2n+1}.$ The representation is given by
$$\lambda (z_1,..,z_{n+1}) = (\lambda^{w_1}z_1,...,\lambda^{w_{n+1}}z_{n+1}).$$
Thus in any factor of $\mathbb C^n=\mathbb C\times...\times\mathbb
C$ the i-th component of the generating field is equal to the
corresponding component of $R$ multiplied by $w_i.$ From this we get
that
$$<V,R>=\sum_{i=1}^{n+1}w_i|z_i|^2 > 0$$ and by Lerman's
theorem we get that $E$ is K-contact. Note that the $S^1-$action on
$S^{2n+1}$ defines a locally free action on $E$ and therefore the
latter is an orbibundle over the symplectic orbifold $\xi
[S^{2n+1}{/}S^1].$ Moreover, for $n=0$ we have the given  Boothby -
Wang fibration, hence a K-contact structure on its total space.

To prove the second part of the theorem, consider separately the
cases $n=1$ and $ n>1.$ In the first case the relevant information
is that the Euler class of the bundle $P\times_{S^1}(\mathbb C
\oplus \mathbb C)$ associated to $P$ by the representation with
weights $w_1,w_2$ is $w_1w_2[\omega ]^2.$ This follows by
calculating the homomorphism induced by the classifying map on
$H^4.$ Then we proceed as in the proof of \ref{prop:example}.

If $n>1,$ then the argument is much simpler, since $M$ has homotopy
3-type of $X,$ thus it also has $b_3$ odd. Again we come to the
contradiction with Theorem \ref{thm:oddbetti} and this completes the
proof.

Note that in any  dimension greater or equal 8 there are  examples
of closed simply connected symplectic manifolds with odd $b_3$ (see,
e.g. \cite{TO},\cite{IRTU}. In particular, the blow up
$\widetilde{\mathbb{C}P}^5$ of $\mathbb{C}P^5$ along the
symplectically embedded Kodaira-Thurston manifold has
$b_3(\widetilde{\mathbb{C}P}^5)=3$ \cite{McD}. More involved are
examples in dimension 8 given by Gompf \cite{G}.

\section{A solution of a problem of Weinstein about fatness}\label{sec:wein}

Let $\Mo$ be a closed symplectic manifold with a Hamiltonian action
of a Lie group $G$ and the moment map $\Psi:M\to \mathfrak g^*$.
Consider the associated Hamiltonian bundle
$$
\Mo \to E:=P\times_G M \to B.
$$

Sternberg \cite{S} constructed a certain closed two--form
$\Omega\in\Omega^2(E)$ associated with the connection $\theta$.  It
is called the {\bf {\em coupling form}} and pulls back to the
symplectic form on each fibre and it is degenerate in general.
However, if the image of the moment map consists of fat vectors then
the coupling form is nondegenerate, hence symplectic.  This was
observed by Weinstein in \cite[Theorem 3.2]{W} where he used this
idea to give a new construction of symplectic manifolds.  In the
sequel, the bundles with a nondegenerate coupling form will be
called {\bf {\em symplectically fat.}} Let us state the result of
Sternberg and Weinstein precisely.

\begin{thm}[Sternberg-Weinstein]
Let $\Mo$ be a a symplectic manifold with a Hamiltonian action of a
Lie group $G$ and the moment map $\nu : M\rightarrow\mathfrak{g}^*$.
Let $G\rightarrow P\rightarrow B$ be a principal bundle. If there
exists a connection in the principal bundle $P$ such that all
vectors in $\nu(M)\subset\mathfrak{g}^*$ are fat, then the coupling
form on the total space of the associated bundle
$$
M\rightarrow P\times_G M\rightarrow B
$$
is symplectic.
\end{thm}
In \cite{W} A. Weinstein asked the following question: {\it are
there symplectically fat fiber bundles whose total spaces  are
simply connected, symplectic but non-K\"ahler?} Note that this
problem was called {\it the Thurston-Weinstein problem for fiber
bundles} in \cite{TO}. The answer to this question is positive and
is a byproduct of  previous sections.

Following Lerman \cite{L3} consider the Kodaira-Thurston manifold
$N$. Recall that the latter is a compact nilmanifold
$(N_3/\Gamma)\times (\mathbb{R}/\mathbb{Z})$, where $N_3$ denotes
the 3-dimensional Heisenberg group of all unipotent $3\times 3$
matrices, and $\Gamma$ is a co-compact lattice in it. It is well
known and easy to see that $N$ is symplectic.

It is known that $N$ can be symplectically embedded into
$\mathbb{C}P^5$ (see, for example \cite{TO}). Consider
$\mathbb{C}P^6$ and define a circle action on it by the formula

$$\lambda\cdot [z_0:z_1:\,...\,:z_6]=
[\lambda z_0:z_1:\,...\,:z_6],\,\lambda\in S^1.\eqno(3)$$
 This
action is trivial on the hyperplane

$$\frak{H}=\{[0: z_1\,...\,:z_6]\,z_i\in\mathbb{C}\}\cong\mathbb{C}P^5.$$
Embed $N$ symplectically into $\frak{H}$. Blow-up $\mathbb{C}P^6$
along $N.$ On the resulting manifold $\widetilde{\mathbb{C}P}^6$ we
have the symplectic circle action extending $(3)$, since  $N$ is
embedded in the fixed point set.

\begin{thm}\label{thm:weinstein} There exist symplectically fat
fiber bundles whose total spaces are simply connected, symplectic
but non-K\"ahler. In particular, the  fiber bundle

$$\widetilde{\mathbb{C}P}^6\rightarrow S^3\times_{S^1}\widetilde{\mathbb{C}P^6}\rightarrow S^2$$
associated with the Hopf bundle $S^1\rightarrow S^3\rightarrow S^2$
admits a fiberwise symplectic structure, but the total space of this
bundle does not admit any K\"ahler structure.
\end{thm}

\noindent{\bf{\em Proof.}} The general theory of symplectic fat
bundles ensures the existence of the non-degenerate coupling form on
the total space. Thus, it is sufficient to prove that the total
space cannot have K\"ahler structures. By a theorem of Lalonde and
McDuff \cite{LMcD}, any hamiltonian fiber bundle over a sphere has
the property of $c$-splitting, that is, the cohomology of the total
space of such bundle is the tensor product (as vector spaces) of the
cohomology of the fiber and the base. It follows that

$$H^*(S^3\times_{S^1}\widetilde{\mathbb{C}P}^6)=H^*(S^2)\otimes_{v.s.}H^*(\widetilde{\mathbb{C}P}^6).$$
It follows from the degree reasons that
$$b_3(S^3\times_{S^1}\widetilde{\mathbb{C}P}^6)=b_3(\widetilde{\mathbb{C}P}^6)=3.$$
Hence the total space cannot be K\"ahler.

\hfill$\square$

\begin{rem}  {\rm Examples of non-simply connected symplectically
fat fiber bundles whose total spaces are not K\"ahler,
were given in \cite{KTW}. According to our knowledge,
simply connected examples of this type were not presented
in the literature, since
techniques like symplectic blow-up construction and the
Lalonde-McDuff $c$-splitting theorem were discovered after the
question was posed. Since the argument presented here is in  the main line of
reasoning in this paper, we decided to write it down.}
\end{rem}

\noindent {\bf Acknowledgement.} The second author started work on the problem described in this paper during his research visit to the Institute of Mathematical Sciences in Chennai. His thanks go to the Institute for wonderful working conditions and especially to Parameswaran Sankaran for valuable discussions. We are also grateful to Jarek K\c edra and Liviu Ornea for answering our questions.


\begin{thebibliography}{ABCDE}

\bibitem[BLG]{BLG} D.E. Blair, S.I.Goldberg, {\it Topology of almost contact
manifolds}, J.Different. Geom. 1(1967), 347--354.


\bibitem[BG]{BG} C.P. Boyer, K. Galicki, {\it Sasakian Geometry},
Oxford Univ. Press, Oxford, 2007.

\bibitem[BGM]{BGM} C.P. Boyer, K. Galicki, B. Mann,
{\it The geometry and topology of 3-Sasakian manifolds}, J. reine angew. Math. 455(1994), 183-220

\bibitem[BMS]{BMS} F. Belgun, A. Moroianu, U. Semmelmann,
{\it Symmetries of contact metric manifolds}, Geom. Dedic.
101(2003), 203-216

\bibitem[BW]{BW} W. M. Boothby, H.C. Wang, {\it On contact manifolds}, Annals Math. 68(1958), 721-734

\bibitem[C]{C} G. Cavalcanti,
{\it The Lefschetz property, formality and blowing-up in
symplectic geometry}, Trans. Amer. Math. Soc. 359(2007), 333-348.

\bibitem[CFM]{CFM} G. Cavalcanti, M. Fern\'andez, V. Mu\~noz,
{\it Symplectic resolutions, Lefschetz property and formality},
Advances Math. 218(2008), 576-599.

\bibitem[F]{F} T. Fujitani, {\it Complex-valued differential forms on
normal contact Riemannian manifolds}, T\^ohoku Math.J. 18(1966),
349--361.


\bibitem[FT]{FT} Y. Felix, J.-C. Thomas, S. Halperin,
{\it Rational homotopy theory}, Springer, Berlin, 2002



\bibitem[GNT]{GNT} O. Goertsches, H. Nozawa, D. T\"oben,
{\it Equivariant cohomology of $K$-contact manifolds}, Math. Ann. 354(2012), 1555-1582.


\bibitem[G]{G} R. Gompf, A new construction of symplectic manifolds,
Annals Math. 142(1995), 527-595

\bibitem[EK]{EK} A. El Kacimi Alaoui, {\it Operateurs transversalement
elliptiques sur un feuilletage riemannien et applications}, Compositio Math. 73(1990), 57-106

\bibitem[K]{K} S. Kobayashi, {\it Principal fibre bundles with the
1-dimensional toroidal group}, T\^ohoku Math. J. 8(1956), 29--45.

\bibitem[KT]{KT} F.Kamber, P.Tondeur,
{\it Duality theorems for foliations}, Asterisque 116 (1984), 458 -
471.

\bibitem[IRTU]{IRTU} R. Ib\'a\~nez, Y. Rudyak, A. Tralle, L. Ugarte,
{\it On certain geometric and homotopy properties of closed symplectic
manifolds},
 Topology Appl. 127(2003), 33-45

\bibitem[KTW]{KTW} J. K\c edra, A. Tralle, A. Woike,
{\it On non-degenerate coupling forms}, J. Geom. Phys. 61(2011),
462-475.

\bibitem[LMcD]{LMcD} F. Lalonde, D. McDuff, {\it Symplectic
structures on fiber bundles}, Topology 42(2003), 309-347

\bibitem[L1]{L1} E. Lerman, {\it How fat is a fat bundle?}  Lett. Math. Phys. 15(1988), 335-339.

\bibitem[L2]{L2} E. Lerman, {\it Contact fiber bundles}, J. Geom. Phys. 49(2004), 52-66

\bibitem[L3]{L3} E. Lerman, {\it A compact symmetric symplectic
non-K\"ahler manifold}, Math. Res. Lett. 3(1996), 587-590

\bibitem[McD]{McD} D. McDuff, {\it Examples of symplectic simply
connected manifolds with no K\"ahler structure}, J. Different. Geom. 20(1984), 267-277

\bibitem[NR]{NR}M. Nicolau, A. Reventos,
{\it Euler class and Gysin sequence for $S^p$-leaves }, Israel J.
Math. 47 (1984), 323--334.

\bibitem[OV]{OV} L. Ornea, M. Verbitsky,
{\it Sasakian structures on $CR$-manifolds}, Geom. Dedicata 125(2007), 159-173

\bibitem[RT]{RT} Y. Rudyak, A. Tralle, {\it On Thom spaces,
Massey products on non-formal simply connected symplectic manifolds},
Internat. Math. Res. Notices, no. 10(2000), 495-513

\bibitem[R]{R} P.Rukimbira, {\it Chern - Hamilton's conjecture and
K-contactness}, Houston J. Math. 21(1995), 709--718.



\bibitem[S]{S} S. Sternberg, {\it Minimal coupling and the
symplectic mechanics of a classical particle in the presence
of a Yang-Mills field}, Proc. Natnl. Acad. Sci. USA 74(1977), 5253-5254


\bibitem[T]{T} S. Tanno, {\it Harmonic forms and Betti numbers of
certain contact Riemannian manifolds}, J. Math. Soc. Japan 19
(1967), 308--316.

\bibitem[TO]{TO} A. Tralle, J. Oprea, {\it Symplectic manifolds with
no K\"ahler structure}, Springer, Berlin, 1997.

\bibitem[WZ]{WZ} Wang Z. Z. and D. Zafran, {\it A remark on the hard
Lefschetz theorem for K\"ahler orbifolds}, Proc. Amer. Math. Soc. 137(2009), 2497-2501

\bibitem[W]{W} A. Weinstein, {\it Fat bundles and symplectic manifolds}, Adv. Math. 37(1980), 239-250


\end{thebibliography}
\end{document}